\input amstex
\documentstyle{amsppt}
\magnification=1000
\hcorrection{0mm} \vcorrection{-5mm}
\hsize=13.2cm \vsize=21.5cm
\pageno=1
\NoBlackBoxes \nopagenumbers \nologo
\headline={ \rightheadline }

\def\blr#1#2{ { \buildrel{#2} \over {#1} } } \def\bur#1#2{ { \underset{#2} \to {#1} }}
\def\blur#1#2#3{{ \blr{\bur{#1}{#2}}{#3} }}

\def\trm#1{ {\text{\rm{#1}}} }
\def\ra{\rightarrow} \def\lra{\longrightarrow} \def\ms{\mapsto} 
   
\def\Pic{{\operatorname{Pic}}}     \def\dim{{\operatorname{dim}}}
     
\def\Hom{{\operatorname{Hom}}}     \def\mult{{\operatorname{mult}}}
\def\un#1{\underline{#1}}

\def\qed{\vskip-3mm \hfill{$\square$} \vskip3mm \parskip=3mm \parindent=4mm}
\def\proof{ \parindent=0mm {{\it Proof}\,:\hskip2mm}\parskip=1mm}
\def\bproc#1{\noindent \bf #1.\hskip2mm\it} \def\eproc{\rm}
\def\rmk#1{{\noindent \bf #1.\hskip2mm}}

\def\lz{\hskip1pt \{ \hskip1pt}          \def\rz{\hskip1pt \} \hskip1pt}
\def\blz{\hskip1pt \big\{ \hskip1pt}     \def\brz{\hskip1pt \big\} \hskip1pt}
\def\bblz{\hskip1pt \bigg\{ \hskip1pt}   \def\bbrz{\hskip1pt \bigg\} \hskip1pt}

\def\pcd{1} \def\dcp{2} \def\foa{3} \def\cnp{4} \def\ncd{5} \def\pnk{6}

\def\tld#1{\widetilde{\hskip1pt#1\hskip1pt}} \def\sg{\vert \sigma \vert}
\def\C{\Bbb C}  \def\Q{\Bbb Q}
\def\CC{{\Cal C}} \def\FF{{\Cal F}} \def\SS{{\Cal S}} \def\NN{{\Cal N}}
\def\pt{\hskip1pt}   \def\"{\lq\lq} \def\eps{\epsilon}
\hyphenation{theorem} \hyphenation{variety}

\def\ptr{\pt \bigstar \pt}    
\def\nts{\pt {}_{{}_\bullet} \pt}

\document
\parskip=3mm \parindent=4mm

 at 8truept
\font\bigfont=cmr12 at 18truept


\ \parskip=1mm \parindent=0mm

\baselineskip 20pt
\centerline{\bigfont An infinitesimal approach}

\centerline{\bigfont to the study of}

\centerline{\bigfont Cycles on Abelian Varieties}

\baselineskip 14pt \vskip4mm

{\bf Giambattista Marini$ \, ^1 $}

{\eightpoint
$ ^1 $\,University of Rome \lq\lq Tor Vergata", via della Ricerca Scientifica, Rome, Italy.
}

\vskip4mm

\centerline{\bf Abstract}

\leftskip15mm \rightskip15mm
\it
This paper is a work in progress on Bloch's conjecture asserting the vanishing of the Pontryagin product of a $ \, p \, $ codimensional cycle on an abelian variety
by $ \, p+1 \, $ zero cycles of degree zero.
We prove an infinitesimal version of the conjecture and we discuss, in particular, the case of $ \, 3 \, $ dimensional cycles.
\rm
\vskip0mm \leftskip0mm \rightskip0mm



{\bf Acknowledgements:}
The author acknowledges the MIUR Excellence Department Project awarded to the Department of Mathematics, University of Rome Tor
Vergata, CUP E83C18000100006.

\vskip3mm
\parskip=3mm \parindent=4mm

In this paper we investigate the geometry of the Chow ring of an abelian variety.
Let $ \, A \, $ be a complex abelian variety of dimension $ \, n \, $ and denote by
$ \, CH_{\bullet}(A) \, $ its Chow group of algebraic cycles modulo rational equivalence,
graded by dimension.
Recall that $ \, CH_{\bullet}(A) \, $ has two ring structures: the first one is given by the intersection product,
which we shall denote by \lq\lq$ \nts $",
the second one is given by the Pontryagin product, which we shall always denote by \lq\lq$ \pt \ptr \pt $".
Observe that the subgroup of zero cycles $ \, CH_0(A) \, $ is a subring of
$ \, CH_{\bullet}(A) \, $ under the Pontryagin ring structure; we shall denote by
$ \, I \, $ the augmentation ideal, namely the ideal of zero cycles of degree zero.
In [Bl], Bloch conjectures a natural generalization to higher codimension of the celebrated square theorem:
if $ \, Y \, \in \, CH_k(A)  $, \, then $ \, I^{\ptr n-k+1} \ptr Y \, = \, 0 $ \
($ I^{\ptr d} \, $ denotes the $ \, d^{th} $-power of $ \, I \, $ under Pontryagin product).
Observe that for $ \, k = n-1 \pt $, \, the cycle $ \, Y \, $ is a divisor and the formula above
particularizes to the cited square theorem:
$ \, Y_{a+b} - Y_a - Y_b + Y \, = \, 0 $, $ \ \forall \ a , \, b \, \in \, A \, $,  \
where $ \, Y_x \, $ denotes the translate of $ \, Y \, $ by $ \, x \, $.
For $ \, k = n \, $ the formula above is trivial.
Bloch proves his conjecture for $ \, k \, = \, 0 , \, 1 , \, 2 \, $. \
After Bloch's results, a great improvement in the understanding of the structure of the Chow ring has been given by Beauville ([Be1], [Be2]):
working with rational coefficients, by the use of Fourier theory he gave a decomposition $ \ CH_\bullet^{}(A) \otimes \Q \, = \, \bigoplus_{d,\,s} [CH_d^{}(A)]_s^{} \ $
and conjectured that negative pieces (those one with $ \, s < 0 $) should be trivial,
where cycles of Beauville's degree $ s $ are defined by putting
\hbox{$ \, [CH_d^{}(A)]_s^{} \, = \, \{ \alpha \, \in \, CH_d^{}(A) \otimes \Q \, \vert \, \mult(m)^* \alpha \, = \, m^{2(n-d)-s} \alpha \pt \} \pt $}.
In spite of the great achievement and clarification about the structure of the Chow ring he obtained, he only proved that negative pieces of he's
decomposition are trivial in the range of dimension $ \, \le 2 \, $ and codimension $ \, \le 1 $, \, the same range of Bloch's results,
outside of this range the main problems are still open. It seems there are no approaches that can manage the higher dimension/codimension situation
(namely dimension $ \, \ge 3 \, $ and codimension $ \, \ge 2 $).
We may stress that the first unknown result concerns 3 dimensional cycles of a 5 dimensional abelian variety.

This paper is a work in progress on Bloch's conjecture:
we prove an infinitesimal version of the conjecture and we discuss, in particular, the case where
$ \, k \, = \, 3 ; $ \ we also show some further results and
connections with the Fourier transform associated with the Poincar\'e bundle.
More in detail, considering the formula
$$
\Pic^0(A)^{\pt k} \, \nts \, CH_k(A) \quad \subseteq \quad I^{\pt \ptr \pt k} \ ,
\tag{\pcd}
$$
we shall prove the following:
\vskip-2mm
$ \bullet \ $ at least for $ \, k = 3 \pt $, \, formula (1) implies Bloch's conjecture \, (see Proposition 16);
\vskip-2mm
$ \bullet \ $ an infinitesimal version of formula (\pcd) holds \, (see Proposition 11).
\vskip-2mm
\noindent
We may note that formula (\pcd) is trivial for $ \, k \le 2 \pt $. \
In fact, the intersection of a divisor in $ \, \Pic^0(A) \, $ with a 1-cycle is a 0-cycle of degree zero
(case $ \pt k = 1 ) ; \ $
furthermore, $ \, I^{\pt \ptr \pt 2} \, $ is the kernel of the abelian sum map $ \, \SS um : \pt I \ra A \, $
and, by rigidity, the intersection between a 2-cycle and two divisors in $ \, \Pic^0(A) \, $
is a zero cycle of degree zero whose abelian sum must be zero
(this take care of the case $ \pt k = 2 ) . $

We also generalize a Mumford's result concerning the countability of components of the rational orbits of the difference map ([Mu]), cfr.\ Lemma (12),
that may help in our situation, cfr. Lemma (14).

\vskip2mm
\noindent
{\S 1 \quad \bf An infinitesimal version of formula (1). }

Let $ \, A \, $ be as above, fix an identification $ \ A \, = \, \C^n / \Lambda \, $ and denote by
$ \, \alpha \pt : \, \C^n \, \ra \, A $ \ the natural projection.
\ Let $ \, W \, $ denote an irreducible and reduced subvariety of dimension $ \, k \, $.
Then consider a very ample linear system $ \ \vert \, D \, \vert $ \ and choose divisors $ \, D_1 , \, ..., \, D_k \, \in \, \vert D \vert \, $ such that
the intersection $ \, D_1 \cap \cdots \cap D_k \cap W \, $ is transverse, in particular it consists of distinct points. Let $ \, r \, $ denote the degree of this intersection, set
$$
\{ p_1 , \, ..., \, p_r \} \quad := \quad D_1 \pt \cap \pt \cdots \pt \cap \pt D_k \pt \cap \pt W
$$
and fix once for all the ordering of the $ \, p_i $'s. \
For a divisor $ \, D \, $ on $ \, A \, $ and $ \, a \in A \, $,
we denote with $ \, t_a D $ \ the translate of $ \, D \, $ by $ \, a \pt $, \, i.e. $ t_a D \, = \, D \ptr \{a\} \pt $.
It is convenient translating by points in the universal covering $ \, \C^n \, $ rather than in $ \, A $,
\, so we also put $ \, t_z D \, = \, t_{\alpha(z)} D \, $ for $ \, z \, \in \, \C^n \pt $.
\ Now let $ \ U \, :=  \, \big\{ z \in \C^n \big\vert \, \vert \vert z \vert \vert \, < \, \eps_0 \, \big\} \ $ and fix $ \, \eps_0 \, $ such that the intersection
$ \, t_{u_1}D_1 \pt \cap \pt \cdots \pt \cap \pt t_{u_k}D_k  \pt \cap \pt W $ \ is transverse for all \ $ u_1, \, ... , \, u_k \, \in \, U $, \,
shrinking $ \, U \, $ if necessary we also assume that the set of all the intersections above consists of $ \, r \, $ connected components.
Each of such components contains exactly one of the $ \, p_i $'s, this fact allow us to define the holomorphic map
$$
\CD
@. U^k @. \ \lra \ @. A^r \\
\un{u} \ := \ @. (u_1^{}, \, ..., \, u_k^{} ) @. \ms @. \big(p_1^{}(\un{u}), \, ..., \, p_r^{}(\un{u}) \big)
\endCD
$$
where the $ \, p_i^{} (\un{u}) \pt $'s are the points of the intersection \ $ t_{u_1^{}}D_1^{} \pt \cap \pt \dots \pt \cap \pt t_{u_k^{}}D_k^{} \pt \cap \pt W $, \,
ordered in such a way to match the previous choice concerning the ordering of the $ \, p_i^{} $'s.

\vskip2mm \noindent
{\bf Local results.}

Fix one of the $ \, p_\iota^{} \, $'s above and denote it by $ \, p \, $. \
As the intersection $ \ D_1^{} \cap \pt \dots \pt \cap D_k^{} \cap W $ \ is $ \, 0 $-dimensional and transverse at $ \, p \pt $,
the tangent space $ \, T_p(A) \, $ decomposes as
$$
T_p(A) \quad = \quad L_1 \oplus ... \oplus L_k \oplus N \ ,
\tag{\dcp}
$$
where
$$
L_j \ := \ \bigcap_{i \neq j} \, T_p(D_i) \cap T_p(W) \qquad \trm{and} \qquad N \ := \ \bigcap_{i} \, T_p(D_i) \ .
$$
We may note that $ \, \dim L_j \pt = \pt 1 \, , \ \forall \, j \ $ and $ \, \dim N \pt = \pt n - k \, $. \  Furthermore,
$$
T_p(D_j) \quad = \quad \bigoplus_{i\neq j} \ L_i \oplus N \quad , \qquad \ T_p(W) \quad = \quad \bigoplus_j \ L_j \ ,
$$
thus $ \, L_j \, $ is a complement of $ \, T_p(D_j) \, $ and $ \, N \, $ is a complement of $ \, T_p(W) \pt $.

Now let $ \, \tilde\pi_j \, $ denote the projections as indicated in the commutative diagram
$$
\CD
T_p(A) @. \ \lra \ L_j \ \subseteq \ @. T_p(A) \\
@V \phi VV @. @VV \phi V \\
\C^n  @. \blr{\lra}{\tilde\pi_j} @. \C^n @. \quad \blr{\lra}{\alpha} \quad @. {\C^n \over \Lambda} \, = \, A \\
\endCD
$$
where the map at the first row is the projection associated with the decomposition (\dcp)
and $ \, \phi \, $ is the natural isomorphism defined by the canonical identifications of vector spaces
\ $ T_p(A) \, = \, T_0(A) \, = \, \C^n \, $.
Similarly, we define $ \, \tilde\pi_N^{} : \, \C^n \ra \C^n \, $ as the morphism induced by the projection
$ \, T_p(A) \ra N \, $. \  Furthermore, we define
$$
\pi_j \ = \ \alpha\circ\tilde\pi_j \quad , \qquad \pi_N^{} \ = \ \alpha\circ\tilde\pi_N^{} \ .
$$
Let $ \, \un{u} \, = \, (u_1, \, ..., \, u_k ) \, \in \, U^k \, $ and $ \, p(\un{u}) \, $ be as previously defined.
Also, consider a series expansion of $ \, p(\un{u}) \, $ and define
$ \, p(\un{u})_i \, := \, i $-order-terms of $ \, p(\un{u}) \, $. \
Then clearly $ \, p(\un{u})_0 \, = \, p \, $ does not depend on $ \, \un{u} \, . $

\bproc{Lemma \foa} We have $ \ p(\un{u})_1 \, = \, \sum_{j=1}^k \pi_j(u_j) \, $,  i.e.
$$
p(\un{u}) \quad = \quad p \ + \ \sum_{j=1}^k \pi_j(u_j) \ + \
\trm{higher-order-terms}
$$
$($in particular, only moving $ \, D_i \, $ in the $ \, L_i $-direction affects the local contribution near $ \, p = p_\iota^{} \, $ of the first order approximation of
the intersection $ \ t_{u_1} D_1 \cap ... \cap t_{u_k} D_k \cap W \pt )$.
\eproc

\proof
As we are interested to the first-order-approximation of $ \, p(\un{u}) \, $,
we are free to work with tangent spaces:
$$
\eqalign{
p(\un{u}) \quad = \quad p \ & + \
\alpha \, \bigg( \big( \phi (T_p(D_1)) + u_1 \big) \, \cap \, ... \, \cap \, \big( \phi (T_p(D_k)) + u_k \big) \, \cap \, \phi (T_p(W)) \bigg) \cr
& + \ \trm{higher-order-terms.}
}
$$
Since by definition $ \ T_p(D_j) \, = \, \oplus_{i\neq j} L_i \oplus N \ $ and $ \ T_p(W) \, = \, \oplus_i L_i \, $,  \  then
$$
\eqalign{
& \tilde\pi_i(u_j) \in \phi\big(T_p(D_j)\big) \, , \ \forall \, i \neq j \, ; \cr
& \tilde\pi_N(u_j) \in \phi\big(T_p(D_j)\big) \, , \ \forall \, j \, .
}
\tag{\foa.1}
$$
By the decomposition $ \, T_p(A) \, = \, L_1 \oplus ... \oplus L_k \oplus N \, $
and by (\foa.1) we obtain
$ \ \phi\big( T_p(D_j) \big) + u_j $ \ = \
$ \phi\big( T_p(D_j) \big) + \tilde\pi_1(u_j) + ... + \tilde\pi_k(u_j) + \tilde\pi_N(u_j) $ \
$ = \ \phi\big( T_p(D_j)\big) + \tilde\pi_j(u_j) \, $.   Therefore,
$$
\eqalign{
& \big( \phi(T_p(D_1)) + u_1 \big) \cap ... \cap \big( \phi(T_p(D_k)) + u_k \big) \cap \phi(T_p(W)) \cr
& = \quad \big( \phi(T_p(D_1)) + \tilde\pi_1(u_1) \big) \cap ... \cap \big( \phi(T_p(D_k)) + \tilde\pi_1(u_k) \big) \cap \phi(T_p(W)) \cr
& = \quad \sum_{j=1}^k \tilde\pi_j(u_j)
}
$$
where the last equality follows again by (\foa.1). This concludes our proof.
\qed


We now introduce some further notation. Denote the $ \, k $-cartesian product of the set $ \, \{0,\, 1\} \, $ by $ \, \Phi \, $,
for $ \ \sigma \, = \, (\sigma_1, \,...,\, \sigma_k) \, \in \, \Phi \ $ and $ \ \un{a} \, = \, (a_1, \, ..., \, a_k) \, \in \, A^k $
\ we define $ \ \un{a}\sigma \, := \, ( ...\, , a_j \sigma_j \, , ... ) \, $,  \  namely
multiplication by $ \, \sigma \, $ replaces $ \, a_j \, $ with $ \, 0 \, $ whenever $ \, \sigma_j \, = \, 0 \, $.  We also define $ \, \sg \, $ as the number of zeros appearing in $ \, \sigma \, $.  As usual, we denote by $ \, Z_0(A) $ \, the free abelian group generated by points of $ \, A \, $
and, for $ \ x \, \in \, A \, $,  \  we denote by $ \, \lz x \rz $ \, the corresponding cycle in $ \, Z_0(A) \, $.  We denote by $ \, \tld{I} \, $ the ideal of $ \, Z_0(A) \, $ of 0-cycles of degree zero and by $ \ {\tld{I}}^{\ptr k} \ $ its $ \, k $-Pontryagin-power
(as far it concerns $ \, 0 $-cycles, the Pontryagin product is defined even if one does not quotient by rational equivalence).
\ Note that $ \ {\tld{I}}^{\ptr k} \ $ is the ideal of $ \, Z_0(A) \, $ generated by the image of $ \, A^k \, $,  \  via the natural map
$ ( a_1 , \, ... , \, a_k ) \, \ms \, \big( \lz a_1 \rz - \lz 0 \rz \big) \ptr ... \ptr \big( \lz a_k \rz - \lz 0 \rz \big) \, . $

By a slight abuse of notation we use the same term \"cycle" to indicate an element in $ \, Z_\bullet(A) \, $ as well as to indicate its associated class in $ \, CH_\bullet(A) \, $.
No confusion will occur since we shall always specify where to consider our \"cycle".

\bproc{Remark \cnp} For $ \, \un{u} \, = \, (u_1, \, ..., \, u_k) \, \in \, U^k \pt $, \
the contribution near $ \, p \, $ of the intersection
$$
\big( t_{u_1}D_1 - D_1 \big) \cap ... \cap \big( t_{u_k}D_k - D_k \big) \, \cap \, W
$$
is just the zero cycle
$$
\zeta_p(\un{u}) \quad := \quad
\sum_{\sigma \in \Phi} \, (-1)^{\sg} \, \blz p(\un{u}\sigma) \brz \quad \in \quad Z_0(A)
$$
\eproc

\noindent
(The word \"near" means \"in the connected component of containing $ p $", \, see the definition of $ \, U \, $ given at the beginning of this section).

\proof
Straightforward.
\qed

It is convenient to our purposes to see our contribution as an element in $ \, \Hom(\Phi,\, A) \pt $, \, namely to factorize $ \, \zeta_p \, $ as
$$
\zeta_p \quad = \quad \delta' \circ \nu_p
\tag{\ncd}
$$
where the maps $ \, \nu_p \, $ and $ \, \delta' \, $ are defined as follows
$$
\CD
U^k @. \quad {\CD @> \nu_p >> \endCD} \quad @. \Hom(\Phi,\, A) @. \quad {\CD @> \delta' >> \endCD} \quad  @. Z_0(A) \\
\un{u} @. \ms @. \big( \sigma \ms p(\un{u}\sigma) \big) \\
@. @. f @. \ms @. \sum_{\sigma \, \in \, \Phi} (-1)^{\sg} \, \blz f(\sigma) \brz
\endCD
$$

The class under rational equivalence of the local contribution $ \, \zeta_p(\un{u}) \, $
might not belong to $ \, I^{\ptr k} \, $,  nonetheless we shall see (proposition (\pnk) below) that the
first-order approximation of $ \ \zeta_p(\un{u}) \ $ even belongs to $ \, {\tld{I}}^{\ptr k} $
\ (as far it concern the approximation to the first order we need not to take the quotient under rational equivalence).

We now come to an elementary proposition, or rather a trivial consequence of \hbox{lemma (\foa)}:

\bproc{Proposition \pnk}
Let $ \ \un{u} \, = \, (u_1, \, ..., \, u_k ) \, \in \, U^k \, $ and let
$ \, p(\un{u})_{\le 1} \, := \, p(\un{u})_0 \, + \, p(\un{u})_1 \, $
be the truncament of $ \, p(\un{u}) \, $ to the first order.  Then
$$
\sum_{\sigma \in \Phi} \ (-1)^{\sg} \,
\blz p(\un{u}\sigma)_{\le 1} \brz \quad \in \quad {\tld{I}}^{\ptr k}
$$
where, as usual, $ \, \ptr \, $ denotes the Pontryagin product.
\eproc

\proof
First of all we may note that by lemma (\foa) and by the definition of the multiplication by $ \, \sigma \, $
the following holds
$$
p(\un{u}\sigma)_1^{} \quad = \quad \sum_{j=1}^k \, \pi_j^{} \big( u_j^{} \sigma_j^{} \big) \quad = \quad
\sum_{j=1}^k \, \pi_j^{} \big( u_j^{} \big) \pt \sigma_j^{}
\tag{\pnk.1}
$$
We have
$$
\eqalign{
\sum_{\sigma \in \Phi} \ (-1)^{\sg} \, \blz p + p(\un{u}\sigma)_1 \brz \quad
& = \quad \lz p \rz \, \ptr \, \sum_{\sigma \in \Phi} \ (-1)^{\sg} \,
\blz p(\un{u}\sigma)_1 \brz \cr
& = \quad \lz p \rz \, \ptr \, \sum_{\sigma \in \Phi} \ (-1)^{\sg} \,
\bblz \sum_{j=1}^k \, \pi_j \big( u_j \big) \pt \sigma_j \bbrz \cr
& = \quad \lz p \rz \, \ptr \, \bigg( \big( \lz\pi_1(u_1)\rz - \lz 0 \rz \big) \, \ptr \, ...
\, \ptr \, \big( \lz\pi_k(u_k)\rz - \lz 0 \rz \big) \bigg) \cr
& \in \quad \lz p \rz \, \ptr \, {\tld{I}}^{\ptr k} \quad = \quad {\tld{I}}^{\ptr k}
}
$$
where the second equality follows by (\pnk.1) and the third equality follows by the straightforward identity
$$
\sum_{\sigma \in \Phi} \, (-1)^{\sg} \, \blz a_j \pt \sigma_j \brz \quad = \quad
\big( \lz a_1 \rz - \lz 0 \rz \big) \, \ptr \, ... \, \ptr \, \big( \lz a_k \rz - \lz 0 \rz \big) \ ,
$$
for all $ \ (a_1, \, ..., \, a_k) \, \in \, A^k \, $. \
\qed

By the previous proof, introducing the map
$$
\CD
\lambda_p \ : \qquad @. T^k @. @>>> @. \Hom(\Phi, \, A ) \\
@. \underline a @. \hskip9mm \ms \hskip-9mm @. @. \Bigg( \dsize \sigma \ \ms \ p + \sum_{j=1}^k \pi_j^{}(a_j^{}) \, \sigma_j^{} \Bigg)
\endCD
$$
where $ T \, = \, \C^n \, $ denotes the universal covering of $ \, A \, $ and $ \, \underline a \, = \, (a_1^{},\,...,\, a_k^{}) $, \, we obtain
$$
\delta' \circ \lambda_p(\underline a) \quad = \quad \{ p \} \, \ptr \, \blur{\ptr}{j\,=\,1}{k} \, \Big( \big\{\pi_j^{}(a_j^{})\big\} - \{0\} \Big) \ .
$$
We want to stress that the various $ \, \pi_j^{} \, $ appearing above depend on decomposition (2), thus they depend on our
$ \, p \, = \, p_\iota^{} \, $.

\vskip2mm
\noindent
{\bf Global results.}

So far we have considered contributions near some point $ \, p \, $ of the intersection at remark (\cnp).
Now we put these contributions all together, namely we define
$$
\eqalign{
\nu \ : \quad & U^k \ {\CD@>\qquad>>\endCD} \ \Hom(\Phi,A)\pt{\vphantom{\big(}}^r \ , \quad
                            \nu(\un{u}) \ := \ \big( \nu_{p_1^{}}^{}(\un{u}) , \, ... , \, \nu_{p_r^{}}^{}(\un{u}) \big) \cr
\zeta \ : \quad & U^k \ {\CD@>\qquad>>\endCD} \ Z_0^{}(A) \ , \quad
                            \zeta(\un{u}) \ := \ \dsize \sum_{\iota=1}^r \ \zeta_{p_i^{}}^{} (\un{u}) \cr
\lambda \ : \quad & T^k \ {\CD@>\qquad>>\endCD} \ \Hom(\Phi,A)\pt{\vphantom{\big(}}^r \ , \quad
                            \lambda(\un{a}) \ := \ \big( \lambda_{p_1^{}}(\un{a}) , \, ... , \, \lambda_{p_r^{}}(\un{a}) \big) \cr
\delta \ : \quad & \Hom(\Phi,A)^r \ {\CD@>\qquad>>\endCD} \ Z_0^{}(A) \ , \quad
                            \delta(f_1^{},\,...,\,f_r^{}) \ := \ \dsize \sum_{\iota=1}^r \ \delta' (f_\iota^{}) \ .
}
$$
We want to stress that by the definition given in remark (\cnp) we have
$$
\zeta(\un{u}) \quad = \quad \big(t_{u_1^{}}^{} D_1^{} - D_1^{} \big) \cap ... \cap \big(t_{u_k^{}}^{} D_k^{} - D_k^{} \big) \cap W
$$
and formula (\ncd) becomes
$$
\zeta \quad = \quad \delta \circ \nu \ ,
\tag{7}
$$
namely we have a commutative diagram
$$
\CD
U^k @> \zeta >> @. \quad Z_0(A) \\
@V \nu VV @. \hskip-16mm \nearrow_{{}_{\ssize \delta}} \\
\Hom(\Phi,A)\pt{\vphantom{\big(}}^r
\endCD
\tag{7$'$}
$$
We may note that $ \, \Hom(\Phi,A)\pt{\vphantom{\big(}}^r \, $ is an abelian variety
(in fact its elements are sets of $ \, r \pt 2^k \, $ points in $ \, A $).
Furthermore, the vertical map is holomorphic and,
in the sense it will be clarified at notation (8) and remark (9) below (that we will introduce also for other purposes) the map $ \, \delta \, $ is Mumford's difference map.

Observe that, by remark (\cnp), the composition of $ \, \zeta \, $ with the natural projection
$ \ cl : \, Z_0(A) \, \ra \, CH_0(A) \ $ is the restriction to $ \, U^k \, $ of the intersection map
$$
\CD
\xi : \quad @. A^k @. {\CD@>\qquad>>\endCD} @. CH_0(A) \\
@. (a_1, \, ...,\, a_k) @. \quad \ms \quad @. \big( t_{a_1}D - D \big) \pt \nts \pt \cdots  \pt \nts \pt \big( t_{a_k}D - D \big) \, \nts \, W
\endCD
$$
(here it is not necessary to distinguish the various $ \, D_j \, \in \, \vert D \vert \pt $, cfr. our former setting)
which is a group-morphism on each factor by the square theorem.

\noindent
{\bf Notation 8. } Let $ \, B \, := \, \Hom(\Phi,A) \pt $, \, we identify
$$
B^r \quad \cong \quad (B^+)^r \times (B^-)^r \quad \cong \quad A^{r \cdot 2^{k-1}} \times A^{r \cdot 2^{k-1}}
$$
where
$ \, B^+ := \, \Hom(\Phi^+,A)$, \ $ B^- := \, \Hom(\Phi^-,A) $,
$ \ \Phi^+ \, := \, \big\{ \, \sigma \in \Phi \, \big\vert \, (-1)^{\sg} = 1 \, \big\} $,
\ $ \Phi^- \, := \, \big\{ \, \sigma \in \Phi \,  \big\vert \, (-1)^{\sg} = -1 \, \big\} \, $
(we may note that $ \, \Phi \, = \, \Phi^+ \cup \Phi^- $).

\noindent
{\bf Remark 9. }
The map $ \, \delta \, $ takes $ \, (f_1^{}, \, ...,\, f_r^{}) \, $ to $ \, \sum_\iota \sum_{\sigma \, \in \, \Phi} (-1)^{\sg} \, \blz f_\iota^{}(\sigma) \brz $,
thus under the previous identifications it is the difference map
$$
\CD
d : \quad @. A^m \times A^m @. \ {\CD @>>> \endCD} \ @. Z_0(A) \\
@. \big( \un{z'} , \, \un{z''} \big) @. \ms @. \sum_{j=1}^m \, \lz z'_j \rz \, - \, \lz z''_j \rz
\endCD
$$
(here $ \, m \, = \, r \pt 2^{k-1} \, $ and the notation is the obvious one: $ \, \un{z'} \, = \, (z'_1,\,...,\,z'_m) \, $ etcetera).


Recalling our identification $ \ A^m \times A^m \ \cong \ \Hom(\Phi^+ \cup \Phi^-,A)\pt{\vphantom{\big(}}^r , \ m \, = \, r \pt 2^{k-1} $, \, let us look at the diagram
$$
\CD
@. \ @. \\
@. U^k @. \ \hookrightarrow \ A^k \ \blr{{\CD@>\qquad>>\endCD}}{\dsize \xi} \ @. @. \ CH_0(A) \ @.
\ \ \blr{{\CD@>\qquad>>\endCD}}{\dsize \mu} \ \ { CH_0^{}(A) \over I^{\ptr k} \phantom{\big(} }
\\
@. @V {\dsize \nu} VV \ \ \ @. @. @AA {\dsize cl} A \\
T^k \quad \blr{{\CD@>\qquad>>\endCD}}{\dsize \lambda} \quad @.
A^m \times A^m @. \blr{{\CD@>\qquad>>\endCD}}{\dsize \delta} @. @. Z_0(A) @. \\
@. \ @.
\endCD
\tag{10}
$$
where $ \ \ \mu \, : \  CH_0(A) \, \lra \, { CH_0(A) \over I^{\ptr k} } \ \ $ is the natural
projection.
In the following proposition we collect a few properties of the diagram (10):

\bproc{Proposition 11} For $ \ a_1^{} , \, ... , \, a_k^{} \, \in \, A \ $ we have
$$
\xi(a_1, \, ...,\, a_k) \quad = \quad \big( t_{a_1}D - D \big) \pt \nts \pt \cdots  \pt \nts \pt \big( t_{a_k}D - D \big) \, \nts \, W \ ,
\tag{11.1}
$$
elements of this type generate $ \, \Pic^0(A)^k \nts CH_k(A) \, $. Furthermore, the following properties hold:

\parskip=1mm \leftskip7mm \parindent-7mm

{\rm a)} \ \ the diagram $(10)$ commutes: $ \ \xi \big\vert_{U^k}^{} \ = \ cl \circ \delta \circ \nu \, ;$

{\rm b)} \ \ $ \xi \, $ is a group-morphism on each factor of $ \, A^k \, = \, A \times ... \times A \, ; $

{\rm c)} \ \ $ \nu \ $ is holomorphic, it is a concrete realization of the intersection at the r.h.s.\ of $(11.1)$,
\ the first order approximation of $ \ \nu(a_1^{} , \, ... , \, a_k^{}) \ $
is $ \ \lambda(a_1^{} , \, ... , \, a_k^{}) $, \ in particular
we have an inclusion of tangent spaces
$ \ \nu_{\ast} \, T_0 \big( U^k \big) \ \subseteq \ \lambda\big(T^k\big) \pt ; $

{\rm d)} \ \ $ \mu \circ cl \circ \delta \circ \lambda \quad = \quad 0 \ $.

\leftskip0mm \parskip=3mm \parindent=4mm
\eproc

\proof
Property (a) has already been proved; property (b) trivially follows by the square theorem;
properties (c) and (d) follow by lemma (\foa) and proposition (\pnk).
\qed

We want to stress that on one hand we have maps
$$
\CD
@. U^k \\
@. @VV \nu V \\
\ T^k @> \lambda >> \ B^r \ @> \ \mu \circ cl \circ \delta \ >> \ { CH_0(A) \over I^{\ptr k} }
\endCD
$$
where the bottom row vanishes, 
on the other hand the composition $ \, \mu \circ cl \circ \delta \circ \nu \, $ is (the restriction of) an additive map on each factor.
As a consequence, by property c),
if $ \, CH_0(A) \, $ had a good topology, or better if its quotient by $ \, I^{\ptr k} \, $ had a good topology, formula (1) would easily follow
(and, in turn, Bloch's conjecture would follow).


We now recall a fundamental result concerning the chow group $ \, CH_0(A) \, $ due to Mumford, \, then we generalize it to the quotient
$ \, CH_0(A) \over I^{\ptr k} \, $.

\bproc{Lemma 12} See {\rm [Mu]}. \ Let $ \, A \, $ be an abelian variety, $ \, m \, $ a positive
integer and consider the difference map (the notation is the usual one)
$$
\CD
d \ : \quad @. A^m \times A^m @. \ {\CD @>>> \endCD} \ @. CH_0(A) \\
@. \big( \un{z'} , \, \un{z''} \big) @. \ms @. \sum_{j=1}^m \, \lz z'_j \rz \, - \, \lz z''_j \rz
\endCD
$$
Then each fibre of $ \, d \, $ is a countable union of algebraic subsets.
\eproc

By observing that $ \, cl \circ \delta \, $ from diagram (10) is the difference map $ \, d \, $ up to the identification
$$
B^r \quad \cong \quad A^{r \cdot 2^{k-1}} \times A^{r \cdot 2^{k-1}}
$$
(cfr. notation 8 and remark 9), \, we obtain the following immediate corollary:

\bproc{Corollary 13}  Each fibre of the composition $ \, cl \circ \delta \, $ from diagram (10) is a countable
union of algebraic subsets.
\eproc

We now come to a generalization of Mumford's lemma.


\bproc{Lemma 14} Let $ \, A \, $ be an abelian variety and let $ \, m , \, k \, $ be integers. Then
each fibre of the composition
$$
A^m \, \times \, A^m \ \ @> \ \ {\dsize d} \ \ >> \ \ CH_0(A) \ \ @> \ \ {\dsize \mu} \ \ >> \ \
{ CH_0(A) \over I^{\ptr k} }
$$
is a countable union of algebraic subsets (c.u.a.s.).
Here, as usual, $ \, d \, $ is the difference map and $ \, \mu \, $ is the natural projection.
\eproc

\proof
The idea of proof is elementary: calling \lq\lq generator" any Pontryagin product of $ \, k \, $ elements in $ \, I $, \,
we show that for $ \, z_0^{} \, \in \, CH_0^{}(A) \, $ and $ \, t \ge 1 \pt $, \, the set of elements in
$ \, A^m \times A^m \, $ whose image in $ \, CH_0^{}(A) \, $ belongs to the same class of $ \, z_0^{} \, $ up
to the sum of $ \, t \, $ {\it generators} \, is a c.u.a.s..

Keeping the notation (8), \, let $ \, z_0 \, \in \, CH_0(A) \pt $, \, let $ \, t \ge 1 \, $ be an integer and consider
$$
\CD
\varphi \ : \ \ @. \ @. \ A^m \ @. \times @. \ A^m \ @. \times @. (A^k)^t @. @. \ \ \lra \ \ @. @.
\ A^m \ @. \times @. {(B^+)}^t @. \times @. \ A^m \ @. \times @. {(B^-)}^t \\
@. \big( @. z_1 @. ; @. z_2 @. ; @. \, ... , \, \un{x_i} , \, ... @. \big) @. \ms @.
\big( @. z_1 @. ; @. \, ... , \, \FF_k^+(\un{x_i}) , \, ... @. ; @.
z_2 @. ; @. \, ... , \,\FF_k^-(\un{x_i}) , \, ... @. \big) @.
\endCD
$$
where $ \, \un{x_1}, \, ..., \, \un{x_t} \, \in \, A^k \, $,
for $ \, \un{x} \, \in \, A^k \, $ one defines $ \, \FF_k^+(\un{x}) \, $ to be the map defined on $ \, \Phi^+ \, $ that takes $ \, \sigma \, $ to $ \, \un{x}\sigma \, $
and $ \, F_k^-(\un{x}) \, $ is defined in the same way (but now the $ \, \sigma $'s run over a different set, namely $ \, \Phi^- \, )$. \
Furthermore, denote by $ \, \pi \, $ the natural projection \ $ A^m \times {(B^+)}^t \times A^m \times {(B^-)}^t \, \lra \, A^m \times A^m  \, $.
Now set
$$
W_t \quad := \quad \big\{ \, ( z_1 , \omega_1; z_2 , \omega_2 ) \, \in \,
A^m \times {(B^+)}^t \times A^m \times {(B^-)}^t \, \big\vert \,
z_1 + \omega_1 \, = \, z_2 + \omega_2 + z_0 \, \in \, CH_0(A) \, \big\}
$$
where $ \, \omega_1 \, $ and $ \, \omega_2 \, $ are considered as effective cycles
under the isomorphisms
$ \, \big(B^+\big)^t \, \cong \, \big(A^{2^{k-1}}\big)^t \, \cong \, \big(B^-\big)^t \, $.
Also set
$$
V \quad := \quad \pi \, \bigg( \, \bigcup_{t=1}^{\infty} \, W_t \cap \trm{image}\, \varphi \bigg) \ .
$$
As $ \, W_t \, $ is a c.u.a.s. by Mumford's lemma, then $ \, V \, $ is a c.u.a.s..
Furthermore, $ \, V \, $ is the fibre $ \, V \, = \, (\mu \circ d)^{-1} ([z_0]) \, $. \
In fact, on one hand the left inclusion \"$\subseteq$" is clear, on the other hand,
if $ \, z_1 - z_2 \, \equiv \, z_0 \ \trm{modulo} \ I^{\ptr k} \, $,  then we can write
$$
\eqalign{
z_1 - z_2 - z_0 \quad & = \quad cl \, \Bigg( \sum_{i=1}^t
\big( \{ a_{i,1} \} - \{ 0 \} \big) \ptr ... \ptr \big( \{ a_{i,k} \} - \{ 0 \} \big) \Bigg) \cr
& = \quad cl \, \Bigg( \sum_{i=1}^t \bigg(
\sum_{\sigma \in \Phi^+} \lz \un{x_i}\sigma \rz \,
- \, \sum_{\sigma \in \Phi^-} \lz \un{x_i}\sigma \rz \bigg) \Bigg)
}
$$
for some $ \, t \, $ and some $ \, \un{x_i} \, \in \, A^k \, $.
\qed

As Mumford's lemma tell us that the fibres of $ \, cl \circ \delta \, $ are c.u.a.s. (corollary 13),
lemma (14) tell us the following:

\bproc{Corollary 15} Each fibre of the composition $ \, \mu \circ cl \circ \delta \, $
from diagram (10) is a countable union of algebraic subsets.
\eproc



\vskip2mm
\noindent
{\S 2 \ \ \ \bf Reduction of Bloch's conjecture to formula (1).}

We now consider Bloch's conjecture for 3 dimensional cycles, namely formula
$$
CH_3(A) \, \ptr \, I^{\ptr n-2} \quad = \quad 0 \ , \qquad \trm{where } \ n \, := \, \dim A \, .
\tag{BC}
$$
We prove that formula (\pcd) implies Bloch's conjecture (BC), \,
namely it is the particular case where $ \, k = r = 0 \, $ of Proposition (16) below.
We remaind that, unfortunately, we were only able to prove that formula (\pcd) holds at the level of first order approximation.

\rmk{Notation}
For $ \, F \, $ and $ \, G \, $ subgroups of $ \, CH_{\bullet}(A) \, $ we denote by $ \, F \nts G \, $ the subgroup generated by intersections
$ \, Y . Z \, $,  where $ \, Y \in F \pt , \ Z \in G \, ; \ $
we also define $ \, F^r \, $ to be $ \ F^r \, := \, F . \cdots . F $ \ $(r $-times).
Furthermore, for a cycle $ \, W \, \in \, CH_{\bullet}(A) \pt $, we let
$ \, W \nts F \, $ (as well as $ \, F \nts W $) denote the group of intersections $ \, W \nts Y \pt , \ Y \in F \pt $.
The groups $ \, F \ptr G \pt , \ F^{\ptr r} \, $ and $ \, W \ptr F \, $ are defined in a similar way.

For the sequel, we fix an ample symmetric divisor $ \, D \pt $, we put $ \, \CC \, = \, D^{n-1} \, $
(namely $ \, \CC \, $ is an ample 1-cycle), and for non-negative integers $ \, j , \, k , \, r \, $ we put
$$
\NN_{j,k,r} \quad := \quad
I^{\ptr j} \, \ptr \, \Big( \Pic^0(A)^k \nts \, CH_3(A) \Big) \, \ptr \, \CC^{\ptr r} \ ,
$$
where $ \, \Pic^0(A)^0 \, = \, [A] \, $ (the unit for the intersection product) \, and
$ \, I^{\ptr 0} \, = \, \CC^{\ptr 0} \, = \, \{o\} \, $ (the unit for the Pontryagin product).

We recall that by Bloch's proposition [Bl, 4.2] we have $ \ \CC \ptr I^{\ptr n} \, = \, 0 \pt $. In fact, more in general:
$$
\CC^{\ptr n-j+1} \ptr I^{\ptr j} = \, 0 \ , \quad j \ = \ 0 \pt , \ 1 \pt , \ 2 \pt , \ ... \pt , \ n+1 \ .
\tag{$ \spadesuit $}
$$

\proof
This result is elementary and well-known, however we remind its proof. As intersection with divisors in $ \, \Pic^0(A) \, $ is a derivation for the Pontryagin product one has
$$
0 \ = \ E_1^{} \nts ... \nts E_j \nts (\CC^{\ptr n+1}) \ = \ {\ssize {n+1 \choose j}} \sigma_1^{} \ptr ... \sigma_j^{} \ptr \CC^{\ptr n+1-j} \ , \quad
$$
where $ \, E_i^{} \, \in \, \Pic^0(A) , \  \sigma_i^{} \, = \, E_i^{} \nts \CC \, \in \, I \pt $.
As $ \, \Pic^0(A) \, $ is divisible, such expressions generate $ \, \CC^{\ptr n-j+1} \ptr I^{\ptr j} \, $ and we are done.
\hfill $ \square $

\bproc{Proposition 16}
Let $ \, A \, $ be an abelian variety of dimension $ \, n \pt $. Assuming
$$
\Pic^0(A)^{\pt 3} \nts CH_3(A) \quad \subseteq \quad I^{\ptr 3}
$$
we then have
$$
\NN_{j,\,k,\,r} \ = \ 0 \ , \ \ \forall \ j + r \ \ge \ n - 2 \qquad ( j , \, k , \, r \ \trm{ non-negative integers}).
\tag{$ \star $}
$$
\eproc

As already mentioned, for $ \, k = r = 0 \, $ (and $ \, j = n-2 $) we obtain Bloch's vanishing
$$
CH_3(A) \ptr I^{\ptr n-2} \, = \, 0 \ .
$$

Defining $ \ \NN_{j,\,k,\,r,\,d}^{} \ $ by replacing $ \, CH_3^{}(A) \, $ with $ \, CH_d^{}(A) \pt $, namely
$$
\NN_{j,k,r,d} \ := \ I^{\ptr j} \ptr \big( \Pic^0(A)^k \nts CH_d(A) \big) \ptr C^{\ptr r} \ = \ 0 \ , \quad
\forall \ j\!+\!r \, \ge \, n\!-\!d\!+\!1 \ ,
$$
proposition (16) is the first not-known case of the following:

\noindent
{\bf Conjecture 17. } Let $ \, A \, $ be an abelian variety of dimension $ \, n $, then formula (1) implies
$$
\NN_{j,k,r,d} \ = \ 0 \ , \qquad \forall \ \ j+r \, \ge \, n-d+1 \ .
$$

It is interesting (after considering rational coefficients) to look at things in the setting of Beauville's theory. First, (1) is equivalent to
the vanishing of
$$
\Pic^0(A)^k \nts [CH_k^{}(A)_\Q^{}]_s^{} \quad = \quad 0 \ , \qquad \trm{for } \ s \le -1 \ ,
$$
vanishing that would trivially follow by Beauville's conjecture stating that
$$
\big[CH_d(A)_\Q^{}\big]_s^{} \quad = \quad 0 \ , \qquad \trm{for } \ s \le -1
$$
($s $ denotes Beauville's degree, cfr.\ [Be1]).
Indeed, Beauville's conjecture would directly imply Bloch conjecture as well the vanishing ($ \star $) from Proposition (16) and, more in general, the vanishing of
$ \, \NN_{j,k,r,d} \, $ (in the correct range, cfr.\ Conjecture 17).

In fact, on one hand, it is not restrictive working with rational coefficients:
as $ \, \Pic^0(A) \, $ is divisible, \, also $ \, I \, $ is divisible, being $ \, I / I^{\ptr 2} \, $ divisible and being $ \, I^{\ptr n+1} \, = \, 0 \pt $.
On the other hand, since $ \ \Pic^0(A)\otimes{\Q} \, = \, [CH_{n-1}(A)_\Q^{}]_1^{} $, \
$ I \otimes {\Q} \, = \, [CH_0(A)_\Q^{}]_{\ge 1}^{} \ $ and $ \ C \otimes 1 \, \in \, [CH_1(A)_\Q^{}]_0^{} $,
our set $ \, {\Cal N}_{j,k,r} \, $ from ($\star$), after tensoring with $ \, \Q \, $, is contained in $ \, [CH_{3-k+r}(A)_\Q^{}]_{\ge j+k} \, $. \
Finally, since $ \, j + r \, \ge \, n - 2 \pt $, \, i.e.\ $ j + k \, > \, n - (3+r-k) \pt $,
assuming that Beauville's conjecture hold we get $ \, [CH_{3-k+r}(A,\Q)]_{\ge j+k} \, = \, 0 \, $. \

We want to stress that Beauville proved his conjecture for $ \, d \, = \, 0,\,1,\,2 $, \, again the case where $ \, d = 3 \, $ is the first case
where the main question is open.


\proof (of prop. 16). \
We shall perform an induction that uses $ \, r \, $ and $ \, k \, $. \
First, we need to deal with the \"boundary" cases. We want to prove the following
\vskip0mm \noindent
$ (16_a) \hskip 30mm (\star) \ $ holds provided that $ \ k \ \ge \ 3 \, $;
\vskip0mm \noindent
$ (16_a) \hskip 30mm (\star) \ $ holds provided that $ \ 3 - k + r \ \ge \ n \, - \, 1 \, $.
\vskip0mm \noindent

Property $ \, (16_a) \, $  holds by trivial reasons of dimension in case $ \, k > 3 \pt $.
In the remaining case, namely the case $ \, k = 3 \pt $,
as by the hypothesis \ $ \Pic^0(A)^3 \nts \, CH_3(A) \, \subseteq \, I^{\ptr 3} $, \, for $ \, j+r \, \ge \, n-2 \ $
we then have $ \, \NN_{j,\,k,\,r} \, \subseteq \, I^{\ptr 3+j} \ptr C^{\ptr r} \, = \, 0 \ $ (by $ \spadesuit $).

As for property $ (16_b) \, $ we first observe that, for trivial reasons of dimension, it holds for $ \, 3 - k + r \, > \, n \pt $.
So we distinguish the two cases $ \, 3 - k + r \, = \, n \, $ and $ \, 3 - k + r \, = \, n - 1 \pt $.
Consider $ \ \sigma \, \in \, {\Cal N}_{j,k,r} \pt $, where $ \, j + r \ \ge \ n - 2 \pt $.
The first case is trivial: as $ \, j + r \, \ge \, n - 2 \, $ and $ \, 3 - k + r \, = \, n \, $,  \  then
$ \ j + k \, \ge \, 1 \, $ and therefore $ \, \sigma \, $ is algebraically trivial.
On the other hand $ \ \dim \, \sigma \, = \, 3 - k + r \, = \, n \, $,  \  so $ \, \sigma \, $ must be trivial.
As for the case where $ \, 3 - k + r \, = \, n - 1 \, $ we first observe that our cycle $ \, \sigma \, $ is a divisor and
that we must have $ \, j + k \, = \, j + 3 + r - n + 1 \, \ge \, 2 \, $ (because $ \, j+r \, \ge \, n -2 $),
then putting together these two observation the thesis follows: we have a divisor
$ \, \sigma \, \in \, I^{\ptr j} \ptr Y $ \ where $ \, Y \, \in \, \big( \Pic^0(A)^k \nts CH_3(A) \big) \ptr \CC^{\ptr r} \, $ and $ \, j + k \, \ge \, 2 \pt $,
\ if either $ \, j \, \ge \, 2 \, $,  or $ \, j \, = \, 1 \, \le \, k \, $ the result follows by the square theorem,
on the other hand, if $ \, k \, \ge \, 2 \pt $, then $ \, Y \, $ is trivial: \, for $ \, D \, $ being an ample divisor,
$ \, R \, \in \, \Pic^0(A)^{k-2} \, $ and $ \, W \, \in \, CH_3(A) \, $, \, the map
$$
\eqalign{
\Psi \, : \ A \, \times \, A \ & \lra \hskip14mm \Pic^0(A) \cr
( \, a \, , \, b \, ) \ & \, \ms \ \big( (D_a - D)\nts(D_b - D)\nts R \nts W \big)
\ptr \CC^{\ptr r} \ ,
}
$$
is trivial on $ \, A \times \{0\} \, \cup \, \{0\} \times A \, $ and therefore it is the zero map.
On the other hand $ \, Y \, $ is clearly a sum of divisors in the image of $ \, \Psi \, $, so it must be trivial.

So far we are done with the \"boundary cases" ($16_a$) and ($16_b$).
\ We are now ready to prove $ \, (\star) \, $.  We proceed by descending induction on $ \, k \, $ and $ \, r \, $,
more precisely we shall prove that
$$
{\Cal N}_{j,k,r} \quad = \quad 0 \ , \qquad \forall \ j + r \, \ge \, n - 2 \, ,
$$
under the hypothesis that $ \, (\star) \, $ holds for $ \, k' \, > \, k \, $ and that $ \, (\star) \, $ hold provided that $ \, k' \, = \, k \, $ and $ \, r' \, > \, r \, $.

By ($16_a$) and ($16_b$) we may assume $ \, k \, \le \, 2 \, $ and $ \, r \, \le \, n - 5 + k \, $. \
In particular, as $ \, j + r \, \ge \, n - 2 \, $ by hypotheses, we obtain $ \, j \, \ge \, 3 - k \, \ge \, 1 \, $.

Let $ \, \zeta \, := \, s_1 \ptr ... \ptr s_j \, $ where $ \, s_i \, = \, F_i \nts \CC \, $ and $ \, F_i \, \in \, \Pic^0(A) \, $,
and let $ \, W \, := \, E_1 \nts \cdots \nts E_{k} \, $,  \  $ E_i \, \in \, \Pic^0(A) \, $. \
Choose $ \, \tld{F_j} \, $ such that $ \, (r+1) \tld{F_j} \, = \, F_j \, $. \  We have
$$
\eqalign{
0 \quad & = \quad
\tld{F_j} \nts \Big( s_1 \ptr ... \ptr s_{j-1} \ptr \big( W \nts Y \big) \ptr \CC^{r+1} \Big) \cr
& = \quad s_1 \ptr ... \ptr s_{j-1} \ptr \Big( \tld{F_j} \nts W \nts Y \Big) \ptr \CC^{r+1} \ + \ \zeta \ptr \big( W \nts Y \big) \ptr \CC^{r} \cr
& = \quad \zeta \ptr \big( W \nts Y \big) \ptr \CC^{r}
}
$$
where the inductive hypotheses take care of the first and the last equality,
and the second equality holds since intersection with a divisor in $ \, \Pic^0(A) \, $ is a derivation for the Pontryagin product.
\qed

\vskip2mm     

\centerline{\bf REFERENCES}

\vskip3mm
\widestnumber\key{MMMM}

\ref
\key Be1
\by A. Beauville
\paper Sur l'anneau de Chow d'une variet\'e ab\'elienne
\jour Math. Ann.
\vol 273
\yr 1986
\pages 647-651
\endref

\ref
\key Be2
\by A. Beauville
\paper Quelques remarques sur la trasformation de Fourie dans l'anneau de Chow d'une variet\'e ab\'elienne
\jour Lecture Notes in Mathematics
\vol 1016
\yr 1983
\pages 238-260
\endref

\ref
\key Bl
\by S. Bloch
\paper Some Elementary Theorems about Algebraic cycles on Abelian Varieties
\jour Inventiones Math.
\vol 37
\yr 1976
\pages 215-228
\endref

\ref
\key Mu
\by D. Mumford
\paper Rational equivalence of zero-cycles on surfaces
\jour J. Math. Kyoto Univ.
\vol 9
\yr 1969
\pages 195-204
\endref

\enddocument